\documentclass[12pt]{amsart}
\usepackage{amssymb,amsmath}
\numberwithin{equation}{section}
\def\dib{\bar\partial}

\def\T{\text}

\def\1#1{\overline{#1}}
\def\2#1{\widetilde{#1}}
\def\3#1{\widehat{#1}}
\def\4#1{\mathbb{#1}}

\def\5#1{\frak{#1}}
\def\6#1{{\mathcal{#1}}}
\def\C{{\4C}}
\def\R{{\4R}}

\def\Z{{\4Z}}

\begin{document}
\abstract
We prove that every  compact, pseudoconvex,  orientable, CR manifold  of $\C^n$,  bounds a complex manifold in the $C^\infty$ sense. 
In particular, $\dib_b$ has closed range.
\endabstract
\title[The range of the tangential Cauchy-Riemann system...]{The range of the tangential Cauchy-Riemann system on a CR embedded manifold}
\author[L.~Baracco]{Luca Baracco}
\address{Dipartimento di Matematica, Universit\`a di Padova, via 
Trieste 63, 35121 Padova, Italy}
\email{baracco@math.unipd.it}
\maketitle
\def\Giialpha{\mathcal G^{i,i\alpha}}
\def\cn{{\C^n}}
\def\cnn{{\C^{n'}}}
\def\ocn{\2{\C^n}}
\def\ocnn{\2{\C^{n'}}}
\def\const{{\rm const}}
\def\rk{{\rm rank\,}}
\def\id{{\sf id}}
\def\aut{{\sf aut}}
\def\Aut{{\sf Aut}}
\def\CR{{\rm CR}}
\def\GL{{\sf GL}}
\def\Re{{\sf Re}\,}
\def\Im{{\sf Im}\,}
\def\codim{{\rm codim}}
\def\crd{\dim_{{\rm CR}}}
\def\crc{{\rm codim_{CR}}}
\def\phi{\varphi}
\def\eps{\varepsilon}
\def\d{\partial}
\def\a{\alpha}
\def\b{\beta}
\def\g{\gamma}
\def\G{\Gamma}
\def\D{\Delta}
\def\Om{\Omega}
\def\k{\kappa}
\def\l{\lambda}
\def\L{\Lambda}
\def\z{{\bar z}}
\def\w{{\bar w}}
\def\Z{{\1Z}}
\def\t{{\tau}}
\def\th{\theta}
\emergencystretch15pt
\frenchspacing
\newtheorem{Thm}{Theorem}[section]
\newtheorem{Cor}[Thm]{Corollary}
\newtheorem{Pro}[Thm]{Proposition}
\newtheorem{Lem}[Thm]{Lemma}
\theoremstyle{definition}\newtheorem{Def}[Thm]{Definition}
\theoremstyle{remark}
\newtheorem{Rem}[Thm]{Remark}
\newtheorem{Exa}[Thm]{Example}
\newtheorem{Exs}[Thm]{Examples}
\def\Label#1{\label{#1}}
\def\bl{\begin{Lem}}
\def\el{\end{Lem}}
\def\bp{\begin{Pro}}
\def\ep{\end{Pro}}
\def\bt{\begin{Thm}}
\def\et{\end{Thm}}
\def\bc{\begin{Cor}}
\def\ec{\end{Cor}}
\def\bd{\begin{Def}}
\def\ed{\end{Def}}
\def\br{\begin{Rem}}
\def\er{\end{Rem}}
\def\be{\begin{Exa}}
\def\ee{\end{Exa}}
\def\bpf{\begin{proof}}
\def\epf{\end{proof}}
\def\ben{\begin{enumerate}}
\def\een{\end{enumerate}}
\def\dotgamma{\Gamma}
\def\dothatgamma{ {\hat\Gamma}}

\def\simto{\overset\sim\to\to}
\def\1alpha{[\frac1\alpha]}
\def\T{\text}
\def\R{{\Bbb R}}
\def\I{{\Bbb I}}
\def\C{{\Bbb C}}
\def\Z{{\Bbb Z}}
\def\Fialpha{{\mathcal F^{i,\alpha}}}
\def\Fiialpha{{\mathcal F^{i,i\alpha}}}
\def\Figamma{{\mathcal F^{i,\gamma}}}
\def\Real{\Re}
%
%
%
\section{  Introduction}
On the boundary  of a relatively compact pseudoconvex domain of $\C^n$, the tangential $\dib_b$ operator has closed range in $L^2$ according to 
Shaw \cite{S85} and Kohn \cite{K86}. 
The natural question arises whether $\dib_b$ has closed range on an embedded, compact CR manifold $M\subset\C^n$ of higher codimension. 
There are two elements in favor of a positive answer. 
 On one hand, by \cite{HL75}, any compact orientable CR manifold of hypersurface type (or maximally complex) $M\subset\C^n$, is the boundary of a complex variety. On the other, by \cite{K86} Section 5, the boundary  of a complex manifold has the property that $\dib_b$ has closed range.
 However, the two arguments do not match: a variety is not a  manifold. A partial answer to the question comes from a different method, that is, the tangential 
H\"ormander-Kohn-Morrey  estimates. They apply to a general, abstract, not necessarily embedded, CR manifold but under the restraint $\dim_{CR}(M)\geq 2$: in this situation, $\dib_b$ has closed range   (Nicoara \cite{N06}). The case of $\dim_{CR}(M)=1$ appears peculiar at first sight: $\dib_b$ does not have closed range in the celebrated Rossi's example of a CR structure in the $3$-dimensional sphere (cf. Burns \cite{B79}).
However, it was conjectured by Kohn and Nicoara in \cite{KN06} that the phenomenon was imputable in full to non-embeddability. We answer in positive to this conjecture  and propose a unified proof of closed range on any embedded CR manifold regardless of its dimension which is solely based upon Kohn's method.
 Precisely, we show that any smooth, compact, pseudoconvex, orientable  CR manifold  embedded in $\C^n$, a boundary in the sense of currents according to Harvey-Lawson, 
is in fact a $C^\infty$ boundary. 
This has an easy explanation in the context of the CR geometry.
 Every such manifold, consists of a single CR orbit (cf. \cite{G69}, \cite{J95}, \cite{MP06}).
Thus, at points of local minimality, which include all points of strong pseudoconvexity, one-sided complexification follows from forced extension of CR functions according to \cite{Tr86}  and \cite{T88}. At points where local minimality fails, this is obtained by propagation along the CR orbit. 
This yields a one-sided complexification of the full $M$, smooth up to the boundary, which is consistent, by  pseudoconvexity, with the portion of the immersed Harvey-Lawson variety which approaches $M$.

I am grateful to Joseph J. Kohn for having given in \cite{K86} the ground of this research, to Emil J. Straube for having attracted  my attention to this problem in its specific approach, and to Alexander E. Tumanov to whose theory of CR minimality my paper is inspired.
\section{ Partial complexification and closed range of $\dib_b$.}
Let $M$ be a smooth, compact manifold of $\C^n$ equipped with the induced CR structure $T^{1,0}M=\C TM\cap T^{1,0}\C^n$. The de-Rham exterior derivative induces  on skew-symmetric antiholomorphic forms a complex that we denote by $\dib_b$. We assume that $M$ is of hypersurface type; thus $\C TM$ is spanned by $T^{1,0}M$, its conjugate $T^{0,1}M$ and a single extra vector field $T$ that we assume to be purely imaginary, that is, satisfying $\bar T=-T$. Let $\gamma$ be a purely imaginary $1$-form which annihilates $T^{1,0}M\oplus T^{0,1}M$ normalized by $\langle \gamma,T\rangle=-1$. $M$ is orientable when there is a global 1-form section $\gamma$ (or vector field $T$). $M$ is pseudoconvex when $d\gamma\geq0$ over $T^{1,0}M\oplus T^{0,1}M$. We will refer to $M$ as ``pseudoconvex-oriented" when it satisfies the combination of the two above properties. 

A CR curve $\gamma$ on $M$ is a real curve such that $T\gamma\subset T^\C M=TM\cap JTM$ where $J$ is the complex structure on $\C^n$. A CR orbit is the union of all piecewise smooth CR curves issued from a point of $M$. According to Sussmann's Theorem (cf. \cite{MP06}) the orbit has the structure of an immersed variety of $\C^n$. Here is a basic, elementary fact
\bp
\Label{p2.1}
(Greenfield, Joricke) 
Let $M\subset\subset \C^n$ be a smooth, compact, connected, CR manifold of hypersurface type. Then $M$ consists of a single CR orbit.
\ep
The result is stated for a hypersurface, the boundary of a domain of $\C^n$; however, its proof readily applies to a CR manifold of hypersurface type (cf. e.g. \cite{MP06} Lemma 4.18). The geometry of the present paper is based upon the following
\bt
\Label{t2.1}
Let $M\subset\subset \C^n$ be a smooth, compact, connected, CR manifold of hypersurface type, pseudoconvex-oriented. Then $M$ is endowed with a partial one-sided complexification in $\C^n$, that is, a complex manifold $X\subset\subset\C^n$ which has $M$ as the smooth connected component of its boundary from the pseudoconvex side.
\et
\br
In the following we will refer to the above circumstance as ``complex extendibility" of $M$ in direction $+JT$ (the positive side being forced by oriented pseudoconvexity). Alternatively, we refer to $X$ as the positive ``partial complexification" in the sense that $TX=TM+\R^+JT$. Notice that the general theory of boundary values of holomorphic functions on a real hypersurface, tells us that if $X$ has a smooth boundary $M$, then $X$ is uniformly smooth up to $M$. This remark underlies all our discussion.
\er
{\it Proof.} \hskip0.2cm
We show that the set of points in whose neighborhood $M$ has a one-sided, positive, partial complexification is non-empty and closed (being trivially open). Take a point $z_o$ of local minimality, that is, a point through which there passes no complex submanifold $S\subset M$; this set of points is certainly non-empty containing, among others, all points of strong pseudoconvexity. 
Take  a local patch $M_o$ at $z_o$  in which the projection $\pi_{z_o}:\C^n\to T_{z_o}M+iT_{z_o}M$ induces a diffeomorphism between $M_o$ and $\pi_{z_o}(M_o)$. Since $\pi_{z_o}(M_o)$ is a  piece of a mininimal  hypersurface, then $(\pi_{z_o}|_{M_o})^{-1}$ extends holomorphically to the pseudoconvex side $\pi_{z_o}(M_o)^+$ by \cite{Tr86} and \cite{T88}, and parametrizes a one-sided complex manifold which has a neighborhood of $z_o$ in $M_o$ as its boundary.
By global pseudoconvexity and by uniqueness of holomorphic functions having the same trace on a real hypersurface, one-sided complex  neighborhoods glue together into a complex neighborhood of a maximal open subset $M_1\subset M$. This is indeed also closed. In fact, let $z_1\in \bar M_1$; since $M$ consists of a single CR orbit by Proposition~\ref{p2.1}, then $z_o$ is connected to any other point $z_o\in M_1$ by a piecewise smooth CR curve $\gamma$. The completion of the proof of the theorem follows from the lemma below.
\bl
\Label{l2.1}
Let $M\subset\subset \C^n$ be a smooth, pseudoconvex-oriented, CR manifold of hypersurface type and let $\gamma$ be a CR curve connecting two points $z_o$ and $z_1$ of $M$. If $M$ has complex extension in direction $+JT(z_o)$ at $z_o$ it has also extension  at $z_1$ in direction  $+JT(z_1)$ .
\el
\bpf
Let $\xi$ be the end-point on $\gamma$ for complex extension and let $\pi_\xi:\C^n\to T_\xi M+iT_\xi M$; then $\pi_\xi(M)$ is a piece of a complex hypersurface
and $\pi_\xi(\gamma)$ is a CR curve. Now, either there is no germ of a complex hypersurface $S$ with $\xi\in S\subset M$ and therefore $(\pi_\xi|_M)^{-1}$ extends holomorphically from $\pi_\xi(M)$ in direction $+\pi'_\xi(JT)$. Otherwise, let such $S$ exist. First, $\pi_\xi(\gamma)$ being a CR curve, it must seat inside $S$ in a neighborhood of $\pi_\xi(\xi)$.
Next, extension of $(\pi_\xi|_M)^{-1}$ to $+\pi'_\xi(JT)$ propagates along $S$ beyond $\pi_\xi(\xi)$ by Hanges-Treves Theorem \cite{HT83}. Thus $\xi=z_1$.

\epf

\bt
\Label{t2.2}
Let $M\subset\subset \C^n$ be a smooth, compact, connected, CR manifold of hypersurface type, pseudoconvex-oriented. Then $\dib_b$ has closed range.
\et
\br
As it has already been said, the theorem is already known when $\dim_{CR}M\geq2$ as a consequence of the tangential H\"ormander-Kohn-Morrey estimates (cf. \cite{N06}). The proof that we give here, in any dimension, is solely based on Kohn method of \cite{K86}.
\er
\bpf
On one hand, $M$ is endowed with  a ``Harvey-Lawson complexification", that is, a   complex, possibly singular,  variety $X$ which has $M$ as boundary in the sense of currents (cf. \cite{HL75}). On the other hand, by identity principle of holomorphic functions, this variety must coincide in a neighborhood of its (immersed) boundary with the non-singular complexification obtained in Theorem~\ref{t2.1}. 
At this stage, the singularities of $X$ are confined to the interior, so they are isolated, and eventually can be removed by desingularization. Altogether, we have obtained a complex manifold $X$ smooth up to the boundary $M$. Thus \cite{K86} Theorem 5.2 can be applied and $\dib_b$ has closed range.

\epf

\end{document}